\documentstyle[amssymb,12pt]{article}
\def\rest{\mathord{\restriction}}
\newcommand{\force}{\Vdash}

\newcommand{\open}{\Bbb}

\newtheorem{theorem}{\sc Theorem}
\newtheorem{proposition}[theorem]{\sc Proposition}
\newtheorem{lemma}[theorem]{\sc Lemma}
\newtheorem{corollary}[theorem]{\sc Corollary}

\newcommand{\hao}{\ha_1}
\newcommand{\goo}{\go_1}
\newcommand{\proof}{\noindent {\sc Proof. }} 
 
\def\mapright#1{\smash{\mathop{\longrightarrow}\limits^{#1}}}

\def\deq{\mathop=\limits^{\rm def}}

\newlength{\labparwidth}
\newcommand{\labpar}[2]{$$\parbox{\labparwidth}{#2}\leqno(#1)$$}

\newcommand{\dom}{\mbox{\rm dom}}
\renewcommand{\hom}{\mbox{\rm Hom}}

\newcommand{\ext}{\mbox{\rm Ext}}

\newcommand{\se}{\subseteq}
\newcommand{\set}[2]{\{#1 \colon #2\}} 
\newcommand{\qed}{$\square$\par\medskip}
\newcommand{\dmd}{\diamondsuit}
\renewcommand{\to}{\rightarrow}
 
\newcommand{\ga}{\alpha}
\newcommand{\gb}{\beta}

\newcommand{\gd}{\delta}

\newcommand{\gz}{\zeta}
\newcommand{\gh}{\eta}

\newcommand{\gk}{\kappa}
\newcommand{\gl}{\lambda}

\newcommand{\gs}{\sigma}
\newcommand{\gt}{\tau}

\newcommand{\gf}{\varphi}

\newcommand{\go}{\omega}
 
\newcommand{\gG}{\Gamma}

\newcommand{\ha}{\aleph}
 
\newcommand{\oP}{{\open P}}
\newcommand{\oQ}{{\open Q}}

\newcommand{\oZ}{{\open Z}}
 
\newcommand{\yd}[1]{y_{\gd, #1}}
\newcommand{\xd}[1]{x_{\eta_\gd(#1)}}
\newcommand{\yt}[1]{y_{\tau, #1}}
\newcommand{\xt}[1]{x_{\eta_\tau(#1)}}
\newcommand{\hh}{\hat{h}}
\newcommand{\kd}[1]{k_{\gd, #1}}
\newcommand{\kt}[1]{k_{\gt, #1}}

\title{Hereditarily Separable Groups and Monochromatic Uniformization}
\author{P. C. Eklof \\University of California, Irvine\thanks{The
 authors thank Rutgers
University for its support.} \and A. H.
Mekler \\Simon Fraser University\thanks{Research partially supported
by NSERC grant \#9848} \and S. Shelah \\Hebrew University
and \\Rutgers University\thanks{Research partially supported by the
BSF.  Publication \#442}}
\date{}

\begin{document}
\maketitle

\begin{abstract}
We give a combinatorial equivalent to the existence of a non-free
hereditarily separable group of cardinality $\hao$. This can be used,
together with a known combinatorial equivalent of the existence of a
non-free Whitehead group, to prove that it is consistent that every
Whitehead group is free but not every hereditarily separable group is
free. We also show that the fact that $\oZ$ is a p.i.d. 
with infinitely
many primes is essential for this result.
\end{abstract}

\section*{Introduction}

An abelian group $G$
 is said to be {\bf separable} if
every finite rank pure subgroup is a free direct summand of $G$;
$G$ is {\bf hereditarily separable} if every subgroup of $G$ is
separable.
It is  well-known that a Whitehead group is hereditarily
separable. In fact, we have the following implications:

\begin{quote}
free $\Rightarrow$  W-group $\Rightarrow$  hereditarily separable 
\end{quote}

Whitehead's Problem asks if the first arrow is reversible. For each of
the two arrows, it has been
proved by the third author that it is independent of ordinary set
theory whether the arrow reverses. (See \cite{Sh44}, \cite{Sh105}
and \cite{Sh98}, or the account in  \cite{EM}.)

Now if we consider the two arrows together, there are four possible
cases, three of which have already been shown to be consistent:

\begin{enumerate}
\item {\em Both arrows reverse.} That is, every hereditarily separable
group is free. This is true in a model of V = L. (See
\cite[VII.4.9]{EM}.)

\item {\em Neither arrow reverses.} That is, there are Whitehead groups
which are not free, and hereditarily separable groups which are not
Whitehead. This is true in a model of MA + $\neg$CH. (See \cite[VII.4.5,
VII.4.6 and XII.1.11]{EM}.)

\item {\em The second arrow reverses but not the first.} That is, every
hereditarily separable group is Whitehead and there are Whitehead
groups which are not free. This is true in a model of Ax(S) +
$\diamondsuit ^*(\omega _1 \setminus S)$ plus $\diamondsuit _\kappa
(E)$ for every regular $\kappa > \aleph _1$ and every stationary
subset $E$ of $\kappa $. (See \cite[Exer. XII.16]{EM}.)

\item {\em The first arrow reverses but not the second.} That, is every
Whitehead group is free, but there are non-free hereditarily separable
groups. It is an application of the main theorem of this paper that this
case is consistent. (See section~\ref{3monu}.)

\end{enumerate}

We also give  additional information about the circumstances under
which Cases 2 and  3 can occur. (See section~\ref{cases 2 and 3}.)
 Finally, we show that Case 4 is impossible for modules over a p.i.d.
with only finitely many (but at least two) primes. (See section~\ref{4monu}.)

 Our methods involve the use of notions of
uniformization, which have played an important role in this subject
since \cite{ShW1}. (See, for example, \cite{EM} and the recent
\cite{EMS}.) It has been proved that there exists a  non-free Whitehead
group of cardinality $\hao$ if and only if there is a ladder system on a
stationary subset of $\goo$ which satisfies $2$-uniformization.
(Definitions are given in detail in the next section.) Our main theorem
here is the following:

\begin{theorem}
\label{main}
A necessary and sufficient condition for
the existence of a non-free hereditarily separable group of cardinality $\hao$
is the existence of a ladder system on a stationary subset of $\goo$ which
satisfies monochromatic uniformization for $\go$ colours.
\end{theorem}

A ladder system $\gh = \{ \gh_\gd \colon \gd \in S \}$
 on $S$ is said to satisfy monochromatic uniformization for $\go$ colours
if for every function $c \colon S \rightarrow \go$, there is a function
$f \colon \goo \rightarrow \go$ such that for every $\gd \in S$,
$f(\gh_\gd(n)) = c(\gd)$ for all but finitely many $n \in \go$.

We believe the main theorem is of independent interest aside from its
use in Case 4.
We will prove sufficiency in section~\ref{1monu} and necessity in
section~\ref{2monu}.
We will then derive the consistency of Case 4 by standard forcing
techniques like those used in \cite{Sh98}. (Actually, we need
only the  sufficiency part of the main theorem for this.) A knowledge of forcing is
required only for sections~\ref{3monu} and \ref{cases 2 and 3}.

We would like to thank Bill Wickless for his help in answering a
question about finite rank torsion-free groups.

\section*{Preliminaries}

We will always be dealing with abelian groups or $\oZ$-modules;  we shall
simply  say ``group''. A  group $G$ is said to be
{\bf $\hao$-free} if  every countable subgroup of $G$ is
free, or equivalently, every finite rank subgroup is free. (See
\cite[IV.2.3]{EM}; throughout the paper we will usually cite
\cite{EM} for results we need, rather than the original source.)
An
$\hao$-free group $G$ is separable if and only if  every pure
 subgroup $H$ of finite rank is a direct summand of $G$, i.e., there is a
projection $h \colon G \rightarrow H$  (a homomorphism such
that $h \rest H$ is the identity on $H$).
The following are two useful facts
 (cf. \cite[IV.2.7 and VII.4.2]{EM}):

\begin{lemma}
\label{Z}
(i) An $\hao$-free group $G$ is separable if every pure cyclic subgroup of
$G$ is a direct summand of $G$.

(ii) An $\hao$-free group $G$ is  hereditarily separable if $B$ is
separable whenever $B$ is a subgroup of $G$ such that $G/B$ is
isomorphic to a subgroup of ${\open Q}/{\open Z}$ and there is a
 finite set $P$ of primes such that the order of every element of
$G/B$
is divisible only by primes in $P$. \qed
\end{lemma}

 A group $G$ is said to be a {\bf Whitehead group} if
$\ext(G, \oZ) = 0$. Every Whitehead group is separable
(\cite[XII.1.3]{EM}).  Since a
subgroup of a Whitehead group is also a Whitehead group, every Whitehead
group is hereditarily separable.

A group   is said to be a {\bf Shelah group} if it has cardinality
$\hao$,  is $\hao$-free, and for every countable subgroup $B$ there is a
countable subgroup $B'\supseteq B$  such that for every countable subgroup
$C$ satisfying  $C \cap  B' = B$, $C/B$ is free.
In that case, we say that $B'$ has the {\bf Shelah property} over $B$.
 In \cite{Sh44} and
\cite{Sh105} it is proved a consequence of Martin's Axiom plus $\neg$CH
that the Whitehead groups of cardinality $\hao$ are precisely the Shelah
groups.

Notions of uniformization (in our sense) were first defined in
\cite{Sh65} and \cite{ShW1}.   Let $S$ be a subset
of $\lim (\omega _1)$. If  $\delta  \in  S$, a {\bf ladder on $\delta
$} is a function $\eta _\delta \colon \omega  \rightarrow \delta $
which is strictly increasing and has range cofinal in $\delta $. A
{\bf ladder system on $S$} is an indexed family  $\eta  = \{\eta
_\delta \colon \delta  \in  S\}$  such that each $\eta _\delta $ is a
ladder on $\delta $.
The ladder system $\gh$ is {\bf tree-like} if whenever $\gh_\gd(n) =
\gh_\gt(m)$, then $ n = m$ and $\gh_\gd(k) =
\gh_\gt(k)$ for all $k < n$.

For a cardinal  $\lambda  \geq  2$, a $\lambda
${\bf -coloring} of a ladder system $\eta $ on $S$ is an indexed family  $c =
\{c_\delta \colon \delta  \in  S\}$  such that  $c_\delta \colon
\omega  \rightarrow  \lambda $.  A {\bf uniformization} of a coloring
$c$ of a ladder system $\eta $ on $S$ is a pair $\langle g,
g^*\rangle $ where
$g\colon \omega _1 \rightarrow  \lambda $, $g^*\colon S \rightarrow
\omega $  and for all  $\delta \in  S$  and all  $n \geq  g^*(\delta
)$, $g(\eta _\delta (n)) = c_\delta (n)$. If such a pair exists, we
say that $c$ can be uniformized. 
We say that $(\eta, \lambda)${\bf -uniformization holds} or that {\bf
$\eta$ satisfies $\gl$-uniformization } if every $\lambda $-coloring
of $\eta $ can be uniformized.

A {\bf monochromatic}
 colouring $c$ of a ladder
system $\eta$ is one such that for each $\gd \in S$, $c_\gd$ is a
constant function. We shall, from now on, consider a monochromatic
colouring with
$\gl$ colours to be a
function $c \colon S \rightarrow \gl$ (which gives the constant value,
$c(\gd)$, of
the colouring of $\gh_\gd$). Then a uniformization of a monochromatic
colouring $c$ is  a pair $\langle f,
f^*\rangle $ where
 for all  $\delta \in  S$  and all  $n \geq  f^*(\delta
)$, $f(\eta _\delta (n)) = c(\delta)$. If every monochromatic
$\gl$-colouring of $\eta$ can be uniformized we say $\eta$ {\bf satisfies
monochromatic uniformization for $\gl$ colours}.

Define a {\bf ladder system based on a countable set} to be an indexed
family $\gh  = \{\gh_\delta \colon \delta  \in  S\}$  such that each
$\gh_\delta $ is a function from $\go$ to a fixed countable set $I$.
We can define notions of colouring and uniformization analogous to those
above. (See \cite[pp. 367--369]{EM}.) The following two results, though
stated and proved for ladder systems on a stationary subset of $\goo$,
are true also for ladder systems based on a countable set.

For ease of reference, we include the following result
 (compare  \cite[XII.3.3]{EM}).

\begin{lemma}
\label{tree-like}

If there is a ladder system on a stationary subset of $\goo$ which
satisfies $\gl$-uniformization (resp. monochromatic uniformization for
$\gl$ colours), then there is a tree-like
ladder system on a stationary subset of $\goo$ which
satisfies $\gl$-uniformization (resp. monochromatic uniformization for
$\gl$ colours).

\end{lemma}

\proof
Suppose $\set{\gh_\gd}{\gd \in S}$ satisfies $\gl$-uniformization (resp. monochromatic uniformization for
$\gl$ colours).
Choose a one-one map $\theta $ from
$^{<\omega }\omega _1$ to $\omega _1$, such that  $\theta (\sigma ) \leq  
\theta (\sigma ')$  if $\sigma '$ is a sequence extending $\sigma $ and such 
that for any  $\tau  \in  {}^{<\omega }\omega _1$,  $\theta (\tau ) \geq
\tau (n)$  for all  $n \in  \dom(\tau )$. Let $C$ be a closed unbounded subset 
of $\omega _1$ consisting of limit ordinals such that for every  $\alpha  \in  
C$, $\theta [^{<\omega }\alpha ] \subseteq  \alpha $. Let  $S' = S \cap  C$.
 For
$\alpha  \in  S'$, define  $\zeta _\alpha (n) = \theta (\langle \eta _\alpha (0)$,
 \ldots
, $\eta _\alpha (n)\rangle)$. Then
 $\{\zeta _\alpha \colon  \alpha  \in  S'\}$ is tree-like and satisfies
 $\gl$-uniformization (resp. monochromatic uniformization for
$\gl$ colours). \qed

\noindent
{\bf Remark. } With a little more care we can prove that if there is a
ladder system on $S$ which
satisfies $\gl$-uniformization, then there is a tree-like
ladder system on the same set $S$ which
satisfies $\gl$-uniformization. (Compare \cite[Exer. XII.17]{EM}.)

\bigskip

The third author has proved that there is  a non-free
Whitehead   group of cardinality
 $\hao$ if and only if there is a
ladder system on a stationary subset of $\goo$ which satisfies
$2$-uniformization.
 (See  \cite[\S XII.3]{EM}.) The main theorem of this paper is an
analogous necessary and sufficient condition for the existence of a
non-free hereditarily separable group. Since every Whitehead group is
hereditarily separable, we can conclude that if there is a
ladder system on a stationary subset of $\goo$ which satisfies
$2$-uniformization, then there is a ladder system on a stationary
 subset of $\goo$ which satisfies
monochromatic uniformization for $\go$ colours. It
is perhaps reassuring to know that there is a simple direct proof of
this consequence:

\begin{proposition}
\label{reassurance}
If there is a
ladder system $\gh$ on a stationary subset of $\goo$ which satisfies
$2$-uniformization, then there is a ladder system on a stationary
subset of $\goo$ which
 satisfies
monochromatic uniformization for $\go$ colours.
\end{proposition}

\proof
By Lemma~\ref{tree-like}, we can assume that $\eta = \{\gh_\gd \colon
\gd \in S \}$ is tree-like.
 Fix a
monochromatic colouring $c$ of $\eta$ with $\go$ colours. Let
$c'$ be a 2-colouring of $\eta$ such that for $\gd \in S$,
and  $k \in \go$ such that $2^k > c(\gd)$,
$$c'_\gd(2^k +1), c'_\gd(2^k +2), \ldots,  c'_\gd(2^{k + 1})$$
 is the sequence $$0^{c(\gd)}, 1, 1, 1, \ldots ,1$$ i.e., $c(\gd)$
 zeroes followed by $2^k - c(\gd)$ ones. Let $\langle g, g^*\rangle$ be
a uniformization of $c'$. To define $\langle f,
f^*\rangle$, let $f^*(\gd)$ be the least $n$ so that $n = 2^{k + 1}$
where $2^k
> \max\{c(\gd)$,
 $g^*(\gd) \}$.
We define $f$
so that for all $\gd \in S$ and $m \geq f^*(\gd)$, $f(\eta_\gd(m)) = c(\gd)$.
To see that $f$ is well-defined, consider the case when
 $\eta_\gt(m) = \eta_\gd(m)$ and $m \geq
f^*(\gt), f^*(\gd)$. Since $\gh$ is tree-like,
 $\eta_\gd(j) =
\eta_\gt(j)$ for all $j \leq m$. By definition of $f^*$, there is a $k$ such that
$2^{k+1} \leq m$, and $2^k
> \max\{c(\gd), c(\gt),
g^*(\gd), g^*(\gt)\}$. But then the values of $g(\gh_\gd(j))$ for $ j =
2^k + 1, \ldots, 2^{k+1}$ code $c(\gd)$ and also $c(\gt)$, so $c(\gd) =
c(\gt)$. \qed

\smallskip

\noindent
{\bf Remark.} The proof actually shows that if $\gh$ is tree-like and
satisfies $2$-uniformization, then $\gh$ satisfies monochromatic
uniformization for $\go$ colours.

 \section{Sufficiency}
\label{1monu}

\begin{theorem}
\label{1monuthm}
If there is a ladder system on a stationary subset  of $\goo$
 which satisfies monochromatic
uniformization for $\go$ colours, then there is a non-free group of
cardinality $\hao$ which is hereditarily separable.
\end{theorem}

\proof
By hypothesis there is a
         stationary subset $S$ of $\go_1$ and a ladder system
$\eta = \set{\eta_\gd}{\gd \in S}$  such that
every monochromatic colouring   with $\go$ colours can be uniformized.
By Lemma~\ref{tree-like}, without loss of generality we can assume that $\eta$ is
tree-like.

We begin by defining the group.
 Let $p_n$ ($n < \go$) be an enumeration of
the primes. The group $G$ will be generated by $\set{x_\ga}{\ga <
\go_1} \cup \set{y_{\gd, n}}{\gd \in S, n < \go}$,
  subject to the
relations 
$$p_n\yd{n+1} = \yd{n} + \xd{n}$$

        For any $\ga$, we let $G_\ga$ denote the subgroup of $G$ generated by
$\set{x_\gb, y_{\gd n}}{\gb < \ga, \gd \in S \cap \ga}$. It is standard
that $G$ is $\ha_1$-free (in fact $\ha_1$-separable) but not free.
        The rest of the proof will be devoted to proving
that
 $G$ is hereditarily
separable. 

Assume $B$ is a subgroup of $G$ such that $G/B$ is
isomorphic to a subgroup of ${\open Q}/{\open Z}$ and there is a
 finite set $P$ of primes such that the order of every element of
$G/B$
is divisible only by primes in $P$;
we need to prove that if  $Z$  is a rank 1
pure subgroup of $B$, then there is a projection of $B$ onto $Z$. (See
Lemma~\ref{Z}.)  Since $G/B$ is
countable there is $\ga$ so that $G_\ga + B = G$ and $Z \se G_{\ga}$.
Fix such an $\ga$ and call it $\ga^*$. Next, choose in $G_{\ga^*}$ a
system of representatives for $G/B$ and let $g\colon G \to G_{\ga^*}$
be the function which assigns to an element  of $G$ its coset
representative. Finally choose $n^*$ so that for all $n \geq n^*$,
$G/B$ is uniquely $p_n$-divisible.

        Let $h_1$ be a projection of $B \cap G_{\ga^*}$ onto $Z$.
Such a projection exists since $B \cap G_{\ga^*}$ is free. We will
extend $h_1$ to a projection $h$ from $B$ to $Z$ by defining $h$ on
$\set{x_\gb - g(x_\gb), \yd{n} -g(\yd{n})}{\gb \geq \ga^*, n \in
\go, \gd \in S, \gd \geq \ga^*}$. Such a definition suffices
(provided it works), since $B$ is generated by this set together with
$B \cap G_{\ga^*}$.

        Define the colouring $c \colon (S \setminus \ga^*) \rightarrow
G_{\ga^*} $ so that $c(\gd)
 = g(\yd{n^*})$. (Since the values are taken in $G_{\ga^*}$, which
is a countable set, this is an allowable colouring.) Let the pair
$\langle f,
f^*\rangle$ uniformize $c$. We can
assume that $f^*(\gd) \geq n^*$ and that $\eta_\gd(f^*(\gd)) \geq
\ga^*$ for all $\gd \in S \setminus \ga^*$. We define the function $h$ in three stages.  First for each
$\gd \geq \ga^*$ and $n \geq f^*(\gd)$, $h(\yd{n} - g(\yd{n})) = 0$ and
$$h(\xd{n} - g(\xd{n})) = h_1(p_ng(\yd{n+1}) - g(\yd{n}) -
g(\xd{n})).$$

        There are two potential problems with the second definition:
namely, why is the right-hand side of the equation defined; and why is
the definition independent of $\gd$? (Note that $\xd{n}$ may equal
$\xt{n}$ for some $\gt$). For the first problem, note
$$p_ng(\yd{n+1}) - g(\yd{n}) - g(\xd{n}) \equiv p_n\yd{n+1} -
\yd{n} - \xd{n} \equiv 0 \pmod{B},$$
hence $p_ng(\yd{n+1}) - g(\yd{n}) - g(\xd{n}) \in B \cap G_{\ga^*}$.
The independence of the definition from $\gd$
is a consequence of the following lemma, after noting that $n \geq f^*(\gd)$
implies $g(\yd{n^*}) = g(\yt{n^*})$ when $\gh_\gd(n) = \gh_\gt(n)$.

\begin{lemma} 
    Let $n \geq n^*$. Suppose that $\tau, \gd \in S$, $\eta_\tau(n) =
\eta_\gd(n)$ and 
$g(\yd{n^*}) = g(\yt{n^*})$. Then $g(\yd{n+1}) = g(\yt{n+1})$.
\end{lemma}

\proof The proof is by induction on $n \geq n^*$.
Since the ladder system is tree-like, we can assume by induction that
$g(y_{\gd, n}) = g(y_{\gt, n})$. Now
$$ p_ng(y_{\gd, n+1})  \equiv  g(y_{\gd, n}) + g(x_{\gh_\gd(n)})
 =  g(y_{\gt, n}) + g(x_{\gh_\gt(n)})
 \equiv  p_ng(y_{\gt, n+1}) \pmod{B}.$$
Since $G/B$ is uniquely $p_n$-divisible (by choice of $n^*$),
$ g(y_{\gd, n+1}) \equiv   g(y_{\gt, n+1}) \pmod{B}$ and hence
$g(y_{\gd, n+1}) =   g(y_{\gt, n+1})$ by definition of
$g$. \qed

        To complete the definition, the second step is to define
$h(x_\gb - g(x_\gb))$ arbitrarily (say 0) for any $\gb$ not covered
in the first step (i.e., $\gb \geq \ga^*$ and $\gb \neq \eta_\gd(n)$
for any $\gd$ and any $n \geq f^*(\gd)$). Finally, for all $\gd$ and $n <
f^*(\gd)$ define $h(\yd{n} - g(\yd{n}))$ as required by the equation
\begin{eqnarray*}
(\yd{n} - g(\yd{n})) + & (\xd{n} - g(\xd{n})) - p_n(\yd{n+1} -
g(\yd{n+1})) + & \\
                       &g(\yd{n}) + g(\xd{n}) -p_ng(\yd{n+1}) &= 0.
\end{eqnarray*}
(Do this by ``downward induction''.)

        It remains to see that $h$ induces a homomorphism. Consider
the free group $F = L \oplus (B \cap G_{\ga^*})$ where $L$ is
the group freely generated by $\set{u_\gb, w_{\gd, n}}{\gd \in S, n\in
\go,
\gb > \ga^*, \eta_\gd(n) > \ga^*}$.
 There is a surjective map $\gf\colon F \to
B$ which is the identity on $B \cap G_{\ga^*}$ and such that $\gf(u_\gb) = x_\gb -
g(x_\gb)$ and $\gf(w_{\gd, n}) = \yd{n} - g(\yd{n})$. The kernel $K$
of $\gf$ is generated by elements of the form $(w_{\gd, n} +
u_{\eta_\gd(n)} - p_nw_{\gd, n+1}) + (g(\yd{n}) + g(\xd{n}) -
p_ng(\yd{n+1}))$. Let $\hh\colon F \to Z$ be defined so that $\hh
\rest B \cap G_{\ga^\ast} = h_1$, $\hh(u_\gb) = h(x_\gb
-g(x_\gb))$ and $\hh(w_{\gd, n}) = h(\yd{n} - g(\yd{n}))$. Since $\hh$
is constantly $0$ on $K$, it induces a homomorphism from $B$ to $Z$
which agrees with $h$ on the generators of $B$.
\qed

\medskip

\noindent
{\bf Remark. } The same proof works with any
tree-like ladder system based on a countable set. (The assumption that
the ladder system is tree-like is necessary, as witnessed by Hausdorff
gaps).
 In
particular, if there is a set of $\ha_1$ branches through the binary tree of
height $\go$ which satisfies   monochromatic
uniformization  for $\go$ colours,
then the group built from these branches is hereditarily separable.
This group is just the group constructed in \cite[VII.4.3]{EM}; there
 it is shown that MA + $\neg$CH implies this group is
hereditarily separable. Given these comments, one might expect that it
is possible to show that MA + $\neg$CH implies that any system of $\ha_1$
branches through the binary tree satisfies monochromatic uniformization
for $\go$ colours.
Indeed, this is the case: given
a set of $\hao$ branches and a monochromatic colouring $c$ by $\go$
colours,
 let the poset, $\oP$, consist of pairs $(s, B)$
where $s$ is a function from ${}^n2 \to \go$ and $B$ is a
finite subset of the branches such that for all $b \in B$, $s(b\rest
n) = c(b)$. If $(t, c) \in \oP$ and
$\dom(t) = {}^m2$,
we define $(t, C) \geq (s, B)$ iff $s \se t$, $B
\se C$ and for all $b \in B$ and $n \leq k \leq m$,
 $t(b\rest k) = s(b\rest n) =  c(b)$. The proof that for each $n$,
$\set{(s, B)}{{}^n2 \se \dom(s)}$ is dense uses the fact that the
colouring is monochromatic. On the other hand the poset is c.c.c.,
since any two conditions with the same first element are compatible.

\smallskip


\section{Necessity}
\label{2monu}

The following lemma can be derived as a consequence of the fact that the
Richman type of a finite rank torsion free group is well-defined (see
\cite{Rich} or \cite {GW}); but for the convenience of the reader we
give a self-contained proof.

\begin{lemma}
\label{finrk}
Suppose $A$ is a torsion free group of rank $r+1$ and every rank $r$
subgroup  is free.  If $B$ and $C$ are pure rank $r$ subgroups,
then the type of $A/B$ is the same as the type of $A/C$.
\end{lemma}

\proof  
The proof is by induction on $r$.  We can assume that $r \geq 1$ and $B
\neq C$.
Consider first the case $r =1$;
then $B \cap C = 0$.
If $b \in B$ and $c \in C$ are generating
elements, then $A \subseteq \oQ b \oplus \oQ c$ and it is enough to
prove that  $m$ divides $b \pmod{C}$ if and only if $m$ divides $c
\pmod{B}$. Now if $m$ divides $b \pmod{C}$, then $b = 
ma + nc$ for some $a \in A$ and $n \in \oZ$. Since $B = 
\langle b \rangle$ is pure in $A$, $m$ and $n$ must be 
relatively prime. Hence there exist $s, t \in \oZ$ such 
that $ns + mt = 1$. But then $c = m(tc - sa) +sb$, so $m$ 
divides $c \pmod{B}$.

Now suppose $r > 1$. Consider $B \cap C$; since $ r + 1 = rk(B +C) = rk(B) +
rk (C) - rk(B \cap C)$ and $2r > r+1$, we have that $rk(B \cap C)
\geq 1$. Since $B\cap C$ is a pure free subgroup of $A$ we can find
$\langle x \rangle \se B\cap C$ which is a pure subgroup of $A$.
Note that $A/\langle x \rangle$ has the property that every subgroup
of rank $r-1$ is free.  Now apply the induction hypothesis to
$A/\langle x \rangle$, $B/\langle x \rangle$ and $C/\langle x
\rangle$. \qed

\begin{theorem}
\label{2monuthm}
If there is a non-free hereditarily separable
 group $G$ of cardinality
$\aleph _1$, then there is a ladder system on a stationary subset of 
$\omega _1$ which satisfies monochromatic uniformization for $\omega $ colours.
\end{theorem}

\proof
We can write $G = \cup _{\alpha < \omega _1} G_\alpha $, a union of a
continuous chain of countable free pure subgroups where, without loss
of generality, we can assume that there is a stationary subset $S$ of
$\omega _1$, consisting of limit ordinals, and an integer $r \geq 0$
such that for all $\delta \in S$, $G_{\delta +1}/G_\delta $ is
non-free of rank $r + 1$ and every subgroup of $G_{\delta
+1}/G_\delta $ of rank $r$ is free. (We use the fact that if a
stationary subset of $\omega _1$ is partitioned into countably many
pieces, then one of the pieces must be stationary: cf.
\cite[II.4.5]{EM}.) Thus for each $\gd \in S$ there is a pure free subgroup $F_\delta
/G_\delta $ of $G_{\delta +1}/G_\delta $ of rank $r$ such that
$M_\delta \deq G_{\delta +1}/F_\delta $ is rank 1 and non-free.
Moreover either $M_\gd $ is divisible or there is a prime $p_\delta $
and an element $y^\delta + F_\delta $ of $M_\delta $ which is not
divisible by $p_\delta $. Without loss of generality (again using
\cite[II.4.5]{EM}), we can assume that there is a prime $p$
such that $p_\delta = p $
 for all $\delta \in S$ such that $M_\gd$ is not divisible.

For each $\gd \in S$, let $\{y^\delta _\ell \colon 0 \leq \ell \leq r\} \se G_{\gd + 1}$ be
 such that $\{y^\delta _\ell + G_\delta \colon 0 \leq \ell \leq r -1\}$
 is a basis of $F_\delta/G_\gd $, $y^\delta _r \notin F_\gd$ and
 $y^\delta _r + F_\delta $ is not divisible by $p$ in $M_\delta $ if
 $M_\gd$ is not divisible. Then $G_{\delta +1}$ is generated by $ G_\gd
 \cup \{y^\delta _\ell \colon \ell \leq r\} \cup \{z^\gd_n\colon n \in
 \omega \}$, where the $z^\gd_n$ satisfy equations $$ r^\delta _nz^\gd_n
 = \sum _{\ell \leq r} s^{\delta ,n}_\ell y^\delta _\ell + g^\gd_n $$

\noindent
where $g^\gd_n \in  G_\delta $ and $r^\delta _n$, $s^{\delta ,n}_\ell
\in  {\open Z}.$ 

Define functions $\varphi _\delta $ on $\omega $ for each $\delta  \in  S$ by:
$$
\varphi _\delta (n) = \langle g^\gd_m,
s^{\delta ,m}_\ell ,  r^\delta _m \colon \ell  \leq  r,  m \leq
n\rangle .
$$

\noindent
Notice that $\varphi _\delta (n)$ determines the isomorphism type of the 
finitely generated subgroup of $G_{\delta +1}/G_\delta $ generated by  (the 
cosets of) $\{y^\delta _\ell \colon \ell  \leq  r\} \cup  \{z^\gd_m\colon m \leq
n\}.$ 

As in \cite[XII.3]{EM} --- see especially Theorem XII.3.3 and the
beginning of the proof of XII.3.1 (p. 381) --- if we show that $\Phi
= \{\varphi _\delta \colon \delta \in S\}$ satisfies monochromatic
uniformization for $\omega $ colours, then there is a ladder system on
a stationary subset of $\omega _1$ with the same property. 

So fix a monochromatic colouring $c\colon S \rightarrow \omega $ of
$\Phi $. We are going to use $c$ to define a subgroup $B$ of $G$ with
a pure cyclic subgroup Z. By the hypothesis on $G$, there will be a
projection $h\colon B \rightarrow Z$.
Because of the way we define $B$ we will be able to use $h$ to define
$f\colon \{\varphi _\delta (n)\colon \delta \in S$, $n \in \omega \}
\rightarrow \omega $ such that for each $\delta \in S$, $f(\varphi
_\delta (n)) = c(\delta )$ for all but finitely many $n \in 
\omega .$  

We will define a continuous chain of subgroups $B_\alpha $ of
$G_\alpha $ by induction on $\alpha $ and let $B = \cup _{\alpha \in
\omega _1} B_\alpha $.  To begin, let $\{x_n \colon n \in \omega \}$
be a basis of $G_0$, and let $B_0$ be the subgroup of $G_0$ generated
by $\{px_0\} \cup \{px_{n+1} - x_n\colon n \in \omega \}$. Thus
$G_0/B_0 \cong Z(p^\infty )$ and $Z \deq {\open Z}px_0$ is a pure
subgroup of $B_0.$ 

Let $A = \{t_n\colon n \in \omega \} \subseteq G_0$ be a complete set
of representatives of $G_0/B_0$ such that $t_0 = 0$.
For each pair
$(d, a)$ where $d > 0$ and $a \in A$, fix
an element $[d, a] \in \go$ such that
 $d t_{[d,a]} + B_0 = a + B_0$.

We will define the $B_\alpha $ so that for all $\ga$

\begin{enumerate}
\item  $B_\alpha  + G_0 = G_\alpha $ and
\item for all $\beta  < \alpha $, $B_\alpha  \cap  G_\beta  = B_\beta .$
\end{enumerate}
Notice then that $G_\alpha /B_\alpha  \cong  G_0/B_0$, and $Z$ is
pure in each  $B_\alpha $.  

The crucial case is when we have defined $B_\delta $ already and
$\delta \in S$. We will define $B_{\delta ,m}$ by induction on $m$
and then let $B_{\delta +1} = \cup _{m \in \omega } B_{\delta, m}$.
Let 
$$ B_{\delta , 0} = \langle B_\delta \cup \{y^\delta _\ell \colon
\ell < r\} \cup \{y^\delta _r - t_{c(\delta )}\}\rangle .  $$

\noindent
Then $B_{\delta ,0} \cap G_\delta = B_\delta $ since $\{y^\delta
_\ell \colon \ell \leq r\}$ is independent mod $G_\delta $. Suppose
$B_{\delta ,m}$ has been defined so that $B_{\delta ,m} \cap G_\delta
= B_\delta $. Thus $(B_{\delta ,m} + G_\delta )/B_{\delta ,m} \cong
G_0/B_0$. Let $ d_m > 0$ be minimal such that $d_mz^\gd_m \in
B_{\delta ,m} + G_\delta $. If $d_mz^\gd_m \equiv a_m \in A$ (mod
$B_{\gd, m}$), let $$B_{\delta ,m+1} = B_{\delta ,m} + {\open
Z}(z^\gd_m - t_{[d_m, a_m]}).$$ Then we will have $B_{\delta ,m+1}
\cap G_\delta = B_\delta $. So, in the end, $B_{\delta +1} \cap
G_\delta = B_\delta $. Moreover, $B_{\delta +1} + G_0 = G_{\delta
+1}$, because, by construction, every generator of $G_{\delta +1}$
belongs to $B_{\delta +1} + G_0$.  

If $\delta \notin S$, the construction of $B_{\delta +1}$ is
essentially the same, except that the colouring $c$ plays no role; we
begin with a set $Y \subseteq G_{\delta +1}$ which is maximal
independent mod $G_\delta $ and let $B_{\delta ,0} = \langle B_\delta
\cup Y\rangle $; then define $B_{\delta ,m}$ by induction as before
(using a well-ordering of type $\omega $ of a set of generators of
$G_{\delta +1}$ mod $G_\delta )$. This completes the description of
the construction of $B.$

Now fix a projection $h\colon B \rightarrow Z$ and fix a
well-ordering, $\prec$, of $Z^{r+1} \times \omega $ of order type
$\omega $. We are going to define the uniformizing function $f$. We
must define $f(\nu )$ for each $\nu $ of the form $\varphi _\delta
(n)$. (Note that there may be many $\delta $ such that $\nu = \varphi
_\delta (n).)$ Suppose 
$$ \nu = \langle g_m, s^m_\ell , r_m \colon
\ell \leq r, m \leq n\rangle .  $$

\noindent
Let $\gs  = \gs (\nu )$  be minimal such that $g_m \in  G_\gs $ for
all $m \leq  n$. 
For each $k \in  \omega $ we can construct a group 
$B^{(k)}_\nu $ just as in the construction of $B_{\delta +1}$, which is 
generated by $B_\gs  \cup  \{y_\ell \colon \ell  < r\} \cup  \{y_r - t_k\}$
together with elements of the form $z_m -
t_{[d_m,a_m]}$ $(m \leq n)$  where the $z_m$ satisfy the relations

\labpar{\ast}{$r_mz_m = \sum _{\ell \leq r} s^m_\ell y_\ell  + g_m$.}

This is an abstract group, which can be regarded as a subgroup of the
free group on $G_\gs  \cup \{y_\ell \colon \ell \leq r\} \cup
\{z_m\colon m \leq n\}$ modulo the relations in $G_\gs$  and the
relations given by $\nu$ (i.e., the equations ($\ast$)). If $\delta $
is such that $\varphi _\delta (n) = \nu $ and $c(\delta ) = k$, then
there is an embedding of $B^{(k)}_\nu $ into $B_{\delta ,n+1}$ which
fixes $B_\gs $ (and is an isomorphism if $\gs = \delta$). As
before, $Z = \langle px_0\rangle $ is a pure subgroup of $B^{(k)}_\nu
$. Since $B^{(k)}_\nu $ is isomorphic to a subgroup of $G$, it is
separable.

Since $h$ exists, there is a  $\prec$-least tuple $\langle w_\ell \colon
\ell \leq r\rangle {\frown} \langle k\rangle $ in $Z^{r+1} \times
\omega $ for which there is a projection $h'\colon B^{(k)}_\nu
\rightarrow Z$ with $h'\rest B_\gs = h\rest B_\gs$, $h'(y_\ell ) =
w_\ell $ for $\ell < r$ and $h'(y_r - t_{k}) = w_r$. Note that this
tuple determines $h'$ on $B^{(k)}_\nu$. Define $f(\nu ) = k.$

We have to show that this definition works, that is, for each $\delta
\in S$, $f(\varphi _\delta (n))$ equals $c(\delta )$ for sufficiently
large $n \in \omega $. Fix $\delta \in S$. With respect to the
well-ordering $\prec$,
there are only finitely
many ``wrong guesses'' which come before the ``right answer'' $\langle
h(y^\delta _\ell ) \colon \ell < r\rangle {\frown} \langle h(y^\delta
_r - t_{c(\delta )})\rangle {\frown} \langle c(\delta )\rangle $. So
we just have to show that no wrong guess can work for all $n$ if it
involves a $k \neq c(\delta )$. If there were a wrong guess that
worked for all $n$ for some $k \neq c(\gd)$, then there would be
a projection $h'$ onto $Z$ whose domain, $B'$, contains $\{y^\delta
_\ell \colon \ell < r\} \cup \{y^\delta _r - t_k\}$, elements of the
form $z^\gd_n - a_{\delta ,n}$ for all $n \in \omega $ (with $a_{\gd, n}
\in A$), and
$B_\sigma $ where $\sigma $ is minimal such that $g^\gd_n \in
G_\sigma$  for all $n \in \omega $.

Let $\tilde{G}$ denote
$$
G_0  + B'  = G_\sigma  +
\langle \{y^\delta _\ell  \colon \ell  \leq  r\} \cup  \{z^\delta_n
 \colon n \in \omega \}\rangle .
$$

\noindent
Notice that for each $g \in G_0$, there is a $j$ such that $p^jg \in B_0
\se
B_\sigma$, which is a subset of $\dom(h)$ and $\dom(h')$. Hence we can
extend $h$ and $h'$ uniquely to homomorphisms from $\tilde{G}$ into
${\open Q}^{(p)} \otimes  Z$. (Here ${\open Q}^{(p)}$ is the
group of rationals whose denominators are powers of
$p.$) Denote the
extension by $\tilde{h}$ (resp. $\tilde{h}'$). We claim that
$\tilde{h} = \tilde{h}'$.

Assume for the moment that this is true. Then
$\tilde{h}(y^{\delta}_r) = \tilde{h}'(y^\delta _r)$. Now
$\tilde{h}(y^\delta _r - t_{c(\delta )}) \in Z$ and
$\tilde{h}'(y^\delta _r - t_k) \in Z$. So
 $$\tilde{h}(t_{c(\delta)} -
t_k) = \tilde{h}(y^\delta _r - t_{c(\delta )}) - \tilde{h}'(y^\delta _r - t_k)
\in  Z.$$
Since $k \neq c(\gd)$, there is an $s \in \oZ
$ so that $s(t_{c(\gd)} - t_k) \equiv x_0
 \pmod{B_0}$, so
 $\tilde{h}(x_0) \in Z$. But $p(\tilde{h}(x_0)) =
\tilde{h}(px_0) = px_0$; this contradicts the fact that $px_0$
generates $Z$.

It remains to prove the claim. Let $H = \tilde{G}/G_\sigma $,
which is isomorphic to $G_{\delta +1}/G_\delta $. Now $\tilde{h} -
\tilde{h}'$ induces a homomorphism from $H$ into ${\open Q}^{(p)}$
since $\tilde{h}$ and $\tilde{h}'$ agree on $B_\gs$, hence on $G_0$
(since $G_0/B_0 \cong Z(p^\infty)$) and so on $G_\gs$.  So it
suffices to prove that $\hom(H, {\open Q}^{(p)}) = 0$.
Assume, to the contrary, that there is a non-zero $\psi
\colon H \rightarrow {\open Q}^{(p)}$.
  Let $K = \ker (\psi )$; then
 the rank of $K$ is
$r$.  Now $H/K$ is isomorphic to a subgroup of ${\open Q}^{(p)}$
and hence is not divisible. By Lemma~\ref{finrk}, $H/K$ is isomorphic
to $M_\gd = G_{\gd + 1}/F_\gd$.  So by the choice of $S$ and $p$,
 $H/K$  is not $p$-divisible;
since $H/K$ is isomorphic to a subgroup of
${\open Q}^{(p)}$, this implies $H/K$ is free. But this is
impossible, since $H$ is not free and $K$ is a subgroup of rank $r$,
and hence free. \qed

\begin{corollary}
If there is an hereditarily separable group of cardinality $\hao$ which
is not free, then there exist $2^{\ha_1}$ different $\ha_1$-separable
groups of cardinality $\hao$ which are hereditarily separable.
\end{corollary}

\proof
By the theorem,  the given hypothesis implies that there is a ladder
system $\gh$ on a stationary subset of $\goo$ which satisfies
monochromatic uniformization for $\go$ colours. Using this ladder
system, we can construct an $\hao$-separable group which is hereditarily
separable as in the proof of Theorem~\ref{1monuthm}. By a standard trick
we can, in fact, construct such groups  with $2^{\hao}$ different
$\gG$-invariants. (Compare \cite[VII.1.5]{EM}.)  \qed

Similarly to the proof of Theorem~\ref{2monuthm} we can prove the
following:

\begin{theorem}
\label{2monuthm bis}

If there is an hereditarily separable
 group $G$ of cardinality
$\aleph _1$ which is not a Shelah group, then there is a ladder system
based on a countable set
 which satisfies monochromatic uniformization for $\omega $ colours.
 \qed
\end{theorem}

\section{Consistency of Case 4}
\label{3monu}

The consistency of Case 4 in the Introduction will now follow from
Theorem~\ref{1monuthm} and the following set-theoretic result.

\begin{theorem}
\label{3monuthm}

It is consistent with ZFC + GCH that the following all hold:

(i) there is a ladder system on a stationary subset of $\goo$ which
satisfies monochromatic uniformization for $\go$ colours;

(ii) there is no ladder system on a stationary subset of $\goo$ which satisfies
$2$-uniformization;

(iii) $\dmd_\gk(E)$ holds for every stationary subset, $E$, of every regular
cardinal $\gk > \hao$.

\end{theorem}


\medskip 
\proof
We assume familiarity with the methods of \cite{Sh98}.
For simplicity  let our ground model be L; fix a stationary, co-stationary subset
$S$ of $\goo$ and a ladder system $\gh$ on $S$. Our forcing $\oP$
will be an iterated forcing with 
countable support using two types of posets: $R$, which adds a Cohen
subset of $\goo$, and $Q(c)$ which is the poset  uniformizing a
monochromatic colouring $c \colon S \rightarrow \go$ of $\gh$ with countable
conditions, i.e.,

$$\begin{array}{lll}
Q(c) & =&\set{f}{f \colon \ga \rightarrow \go \hbox{ for some
successor } \ga <
\goo \hbox{ and for all } \gd \in S \cap \ga,\\ 
&& f(\gh_\gd(n)) = c(\gd) \hbox{ for almost all } n \in \go}.
\end{array}$$

In the iteration $\oP$ we force with $\tilde{R}$ at successors of
even ordinal stages and force with $\tilde{Q}(\tilde{c})$ at the
successors of odd ordinal stages, where, as usual, the
names $\tilde{c}$ are chosen so that all possibilities occur. The
posets $R$ and $Q(c)$ are proper, so stationary sets are preserved by
$\oP$. Also, $\oP$ is $(\goo \setminus S)$-closed and of cardinality
$\ha_2$, so GCH holds in the generic extension as well as
$\dmd_\gk(E)$ for every stationary subset of every regular cardinal
$\gk > \hao$.

It remains to show that in the generic extension $2$-uniformization
fails  for every stationary subset $E$ of $\goo$ and every ladder system
$\gz = \set{\gz_\gd}{\gd \in E}$. By doing an initial segment of the
forcing we can assume that $E$ and $\gz$ are both in the ground model.
Let $X$ be the generic set for the first copy of $R$ in the iteration of
$\oP$. Consider the $2$-colouring $\set{c_\gd}{\gd \in E}$ of $\gz$
defined as follows: $c_\gd(n) = 0$ if and only if $\gd + n \in X$.  The
proof that this colouring is not uniformized now follows along the same
lines as that in \cite{Sh98}. \qed

\begin{corollary}
\label{3monucor}

It is consistent with ZFC + GCH that there is an hereditarily separable
group of cardinality $\hao$ which is not free, and every Whitehead group
(of arbitrary cardinality)
is free.

\end{corollary}

\proof
We use the model of ZFC + GCH constructed in Theorem~\ref{3monuthm}.
        Clause (i) in Theorem~\ref{3monuthm} together with
 Theorem~\ref{1monuthm} imply that
there is a non-free hereditarily separable group of cardinality
$\hao$.  Clause (ii) implies that there is no non-free Whitehead
group of cardinality $\hao$.
(See \cite[XII.3.1(i)]{EM}.)
 Finally clause (iii) enables one to do
an inductive proof that there is no non-free Whitehead group of any
cardinality (as, for example, in \cite[XII.1.6]{EM}). \qed


%

\section{Cases 2 and 3}
\label{cases 2 and 3}

In Cases 2 and 3 of the Introduction we are  in the situation where there is
a Whitehead group which is not free; here we shall consider  two
 hypotheses which are stronger than this hypothesis:
 first, that there is a Whitehead group which is
not a Shelah group; and, second, that every Shelah group is a Whitehead
group.

The following theorem says that the hypothesis that there
is a Whitehead group which is not a Shelah group is not
consistent with Case 3. It also gives another consistency proof for Case 2
 since it is known that it is consistent that there are
Whitehead groups which are not Shelah groups (see \cite[XII.3.11]{EM}).

\begin{theorem}
\label{case3}
If there is a Whitehead group of cardinality $\hao$ which is not a Shelah
group, then there is an hereditarily separable group which is not a
Whitehead group.
\end{theorem}

\proof
By \cite[XII.3.19]{EM} there is a ladder system $\gh =
\set{\gh_\gd}{\gd \in \goo}$ based on a countable set
$I$ which satisfies 2-uniformization. Without loss of generality we can
assume that $I = \go$ and each $\gh_\gd \colon \go \to \go$ is strictly
increasing. Moreover, as in the proof of Lemma~\ref{tree-like}, we can
assume that $\gh$ is tree-like, and hence, as in the proof of
Lemma~\ref{reassurance}, $\gh$ satisfies monochromatic uniformization
for $\go$ colours.

For each $\gd \in \goo$ and $n \in \go$ let $$k_{\gd, n} = (\gh_\gd(n) +
1)!$$  and
$$ k_{\gd, n}' = \gh_\gd(n)! $$
Let $G$ be the group generated by $\set{x_n}{n <
\go} \cup \set{y_{\gd, n}}{\gd \in \goo, n < \go}$,
  subject to the
relations 
\begin{equation}
\kd{n+1}\yd{n+1} = \yd{n} + \xd{n}. \label{g}
\end{equation}
As in the proof of Theorem~\ref{1monuthm}, $G$ is hereditarily
separable.

 It remains to show that $G$ is not a Whitehead group. For this we shall
 define an epimorphism $\pi \colon H \to G$ with kernel $\oZ$ which does
 not split. Let $H$ be the group  generated by
 $\set{x_n'}{n <
\go} \cup \set{y_{\gd, n}'}{\gd \in \goo, n < \go} \cup \{z\}$,
  subject to the
relations 
\begin{equation}
\kd{n+1}\yd{n+1}' = \yd{n}' + \xd{n}' + \kd{n+1}'z. \label{h}
\end{equation}
There is an epimorphism $\pi$ taking $\yd{n}'$ to $\yd{n}$, $x_m'$ to
$x_m$, and $z$ to $0$; the kernel of $\pi$ is the pure subgroup of $H$
generated by $z$. Aiming for a contradiction, assume there is a
splitting of $\pi$, i.e., a homomorphism $\gf \colon G \to H$ such that
$\pi \circ \gf = 1_G$.
So $\gf(\yd{n}) - \yd{n}' \in \ker(\pi)$ for all $\gd < \goo$, $n \in
\go$.
 Since a countable union of countable sets is
countable, there exists $\gd \neq \gt$ such that $\gh_\gd(0) =
\gh_\gt(0)$, $\gh_\gd(1) =
\gh_\gt(1)$ and $\gf(\yd{0}) - \yd{0}' = \gf(\yt{0}) - \yt{0}'$. Let $m$
($ \geq 2$) be minimal such that $\gh_\gd(m) \neq \gh_\gt(m)$.

We claim that $\gf(\yd{n}) - \yd{n}' = \gf(\yt{n}) - \yt{n}'$ if $n <
m$. The proof is by induction on $n < m $; the initial case $n = 0$ is
by choice of $\gd$ and $\gt$. So supposing the result is true for
 $n < m - 1$, we will prove it for $n + 1$. Applying the homomorphism
 $\gf$ to equation (\ref{g}) for $\gt$ as well as $\gd$ and subtracting we
 get that (in $H$)
 \begin{equation}
 \label{i}
 \kd{n+1}\gf(\yd{n+1}) - \kt{n+1}\gf(\yt{n+1}) = \gf(\yd{n}) -
 \gf(\yt{n})
 \end{equation}
 since $\xd{n} = \xt{n}$ because $n < m$. But then by induction
\begin{equation}
 \label{ii}
 \kd{n+1}\gf(\yd{n+1}) - \kt{n+1}\gf(\yt{n+1}) = \yd{n}' -
\yt{n}'.
 \end{equation}
Now by equation (\ref{h}),  since $\xd{n} = \xt{n}$
and $\kd{n+1}' = \kt{n+1}'$ (the latter because $n < m - 1$), we have
\begin{equation}
\label{iii}
 \kd{n+1}\yd{n+1}' - \kt{n+1}\yt{n+1}' = \yd{n}' -
\yt{n}'.
 \end{equation}
so by equations (\ref{ii}) and (\ref{iii}) we have
\begin{equation}
\label{iv}
\kd{n+1}(\gf(\yd{n+1}) - \yd{n+1}') =
\kt{n+1}(\gf(\yt{n+1}) - \yt{n+1}').
\end{equation}
Since $n < m - 1$, $\kd{n+1} = \kt{n+1}$, so cancelling $\kd{n+1}$ from
equation (\ref{iv}), we obtain the desired result, and the claim is
proved.

Now equation (\ref{ii}) holds for $n = m - 1$ so
\begin{equation}
\label{v}
 \kd{m}\gf(\yd{m}) - \kt{m}\gf(\yt{m}) = \yd{m-1}' -
\yt{m-1}'.
 \end{equation}
In this case, instead of (\ref{iii}) we have
\begin{equation}
\label{vi}
 \kd{m}\yd{m}' - \kt{m}\yt{m}' - (\kd{m}' - \kt{m}')z = \yd{m-1}' -
\yt{m-1}' .
 \end{equation}
so combining (\ref{v}) and (\ref{vi}) we have
\begin{equation}
\label{vii}
 \kd{m}\gf(\yd{m}) - \kt{m}\gf(\yt{m}) =
  \kd{m}\yd{m}' - \kt{m}\yt{m}' - (\kd{m}' - \kt{m}')z.
\end{equation}
Say $\gh_\gd(m) < \gh_\gt(m)$. Then $\kd{m}$, $\kt{m}'$ and $\kt{m}$ are
all divisible by $\kd{m} = (\gh_\gd(m) + 1)!$
so equation (\ref{vii}) implies that $(\gh_\gd(m) + 1)!$ divides $\kd{m}'z =
\gh_\gd(m)!z$ in $H$ which is a contradiction, since $z$ generates a
pure subgroup of $H$. \qed

Now we consider the hypothesis that
every Shelah group  is a Whitehead group.
This is true in a model of Martin's Axiom, in which case
there are hereditarily separable groups which are not Whitehead groups,
i.e., Case 2 holds.
Here we show that it is consistent that every Shelah group
is a Whitehead group but every hereditarily separable group is a
Whitehead group, i.e., there is a model for Case 3 in which every Shelah
group is a Whitehead group. For this purpose
we use the notion of  stable forcing. A poset, $\oP$, is
{\bf stable} if for every countable subset $P_0$ there is a
countable subset $P_1$ so that for every  $p \in \oP$ there is an
extension  $p'$ of $p$  and an element  $p^* \in P_1$ so that $p'$ and $p^*$
are compatible with exactly the same elements of $P_0$. In \cite{AS} the basic
facts about c.c.c.\ stable forcings are proved.
There are a few basic facts that we will use:

\begin{proposition}
\label{stable}
\begin{enumerate}
\item \cite{AS} Any iteration of c.c.c.\ stable forcings with finite
support is c.c.c.\ and stable.

\item The forcing adding any number of Cohen reals is   stable.

\item If $A$ is  a Shelah group and
$$0 \rightarrow  \oZ \mapright{} B \mapright{\pi } A \rightarrow
0$$
is a short exact sequence, then the finite forcing, $Q(\pi)$,
constructing the
 splitting of $\pi$  is (c.c.c. and) stable. \end{enumerate} \end{proposition}

\proof We will prove only  the last of the statements. Write $A$ as
$\bigcup_{\ga < \go} A_\ga$ (an $\goo$-filtration) where
each $A_\ga$ is pure in $A$ and $A_{\ga +1}$
has the Shelah property over $A_\ga$. The forcing $Q(\pi)$ is the set of
partial splittings of $\pi$ whose domains are finite rank pure subgroups of $A$.
(This forcing is c.c.c. --- see, e.g., \cite[XII.1.11]{EM}.) Given
$P_0$,
 choose  $\ga$ so that every element of $P_0$ has
domain contained in $A_\ga$. Let $P_1$ be the set of
elements of $Q(\pi)$ whose domains are contained in $A_{\ga+\go}$.
Given $p \in Q(\pi)$, let $G$ be the pure subgroup of
$A$ generated by $A_\ga \cup \dom(p)$. There exists $n \in \go$ such that $G
\cap A_{\ga + \go} = G \cap A_{\ga + n}$ (since $G$ has finite rank over
$A_\ga$). Then $$G = (G \cap A_{\ga +n}) \oplus \langle y_0, \ldots,
y_m\rangle$$ for some $y_0, \ldots, y_m$ since $G \cap A_{\ga + n + 1} =
G \cap A_{\ga + n}$ and $A_{\ga + n + 1}$ has the Shelah property over
$A_{\ga + n}$. Extend $p$ to $p' \in Q(\pi)$
such that
 $$\dom(p') = M \oplus \langle y_0, \ldots,
y_m\rangle$$
where $M \se G \cap A_{\ga + n}$ 
is a finite rank pure subgroup of $G \cap A_{\ga + n}$ 
such that $\dom(p) \se M \oplus \langle y_0, \ldots,
y_m\rangle$. Let $p^* = p'\rest M \in P_1$.

It suffices to prove that if $q \in P_0$ is compatible with $p^*$, then
$q$ is compatible with $p'$. So suppose that $r \in Q(\pi)$ such that $r
\geq q, p^*$. Without loss of generality $\dom(r) \se G \cap A_{\ga +
n}$. Define $r'$ with $$\dom(r') = \dom(r) \oplus \langle y_0, \ldots,
y_m\rangle$$ by: $r'\rest \dom(r) = r$ and $r' \rest \langle y_0, \ldots,
y_m\rangle = p' \rest \langle y_0, \ldots,
y_m\rangle$. Clearly $r' \geq q, p'$. Moreover, $\dom(r')$ is pure in
$A$ since $G$ is pure in $A$ and $\dom(r)$ is pure in
$G \cap A_{\ga + n}$; so $r'\in Q(\pi)$. \qed
\begin{theorem}
\label{case3bis}
It is consistent that every Shelah group is a Whitehead group and every
hereditarily separable group is a Whitehead group.
\end{theorem}

\proof We do our forcing over L by iteratively adding  subsets of
$\goo$ by finite conditions and adding splittings for Shelah groups.
More precisely, our forcing $\oP$
will be an iterated forcing with 
finite support and of length $\go_2$ using two types of posets: $R$,
the finite functions from $\goo$ to $2$, and $Q(\pi)$ which is the
finite forcing splitting  $\pi$ as in Proposition~\ref{stable}(3).
 If we choose
the iterants correctly, then in the generic extension
every Shelah group of cardinality $\hao$ will be a Whitehead group and
$\dmd(E)$ will hold for every stationary subset of every regular cardinal
greater than $\hao$. It will suffice then to show that every hereditarily separable
group of cardinality $\hao$ is a Shelah group (because we have all
instances of diamond above $\hao$: cf. \cite[Exer. XII.16(ii)]{EM}.)

By Theorem~\ref{2monuthm bis} it is enough to show that, in the generic
extension, if $\Phi =
\{\gf_\ga\colon \ga < \goo\}$ is a ladder system based on $\go$,
 then  $\Phi$ does not satisfy monochrome
uniformization for $\go$ colours. In fact, we will show that $\Phi$ does
not satisfy monochrome uniformization for $2$ colours.  By absorbing an
initial segment of the forcing into the ground model we can assume that
$\Phi$ is in the ground model and the
forcing $\oP$ is first $R$, the finite functions from $\goo$ to $2$, followed
by a name $T$ for a c.c.c.\ stable forcing.
 We define the colouring $c \colon \goo \to 2$
to be  the generic set for $R$;
 let $\tilde{c}$ be a name for $c$.

In order to obtain a contradiction, assume that this colouring can be
uniformized. Then there is a  pair $\langle \tilde{f}, \tilde{f}^*\rangle$
of names for functions and there is a $p'
\in \oP$  such that $p' \force ``\langle \tilde{f}, \tilde{f}^*\rangle
\hbox{ uniformizes } \tilde{c}$''.
 Now let $P_0$ be a countable
subset of $\oP$ containing $p'$
 as well as
 for every $n < \go$,
 a maximal
antichain which determines the value of $\tilde{f}(n)$. Let $P_1$ be as given
by the definition of a stable poset for this $P_0$. For each $\ga \in
\goo$
choose $p_\ga \geq p'$ so
that $p_\ga$ determines the values of $\tilde{f}^*(\ga)$ and
$\tilde{c}(\ga)$
 and there exists $p^*_\ga
\in P_1$ so that $p_\ga$ and $p^*_\ga$ are compatible with exactly
the same elements of $P_0$.
Say $$p_\ga \force   \tilde{f}^*(\ga) =
m_\ga \wedge \tilde{c}(\ga) = e_\ga.$$
 By the pigeon-hole
principle, there exists an uncountable set $E\se \goo$ and $p^*\in P_1$
so that for all $\ga \in E$, $p_\ga^* = p^*$.
Since $p^*$ is compatible with $p'$, there exists $q_1 \geq p^*, p'$.
By the definition of $R$, there exists $\ga_0 \in E$ and $q_2 \in \oP$
such that $q_1 \leq q_2$ and $q_2 \force \tilde{c}(\ga_0) \neq
e_{\ga_0}$. So $q_2 \force ``\exists k > m_{\ga_0} \hbox{ s.t. }
\tilde{f}(\gf_\ga(k)) \neq e_{\ga_0}$''.
Thus there exists $q_3 \geq q_2$ and $k_0 > m_{\ga_0}$ such that $q_3
\force \tilde{f}(\gf_\ga(k_0)) \neq e_{\ga_0}$. But there is a maximal
antichain in $P_0$ of conditions forcing the value of
$\tilde{f}(\gf_\ga(k_0))$. Hence there exists $r \in P_0$ and $q_4$ such
that $r \leq q_4$, $q_3 \leq q_4$ and $r \force \tilde{f}(\gf_\ga(k_0)) =
1 - e_{\ga_0}$. Then $r$ is compatible with $p^* = p^*_{\ga_0}$ and
hence with $p_{\ga_0}$. But this
is a contradiction since $p_{\ga_0} \force \tilde{c}(\ga) = e_{\ga_0}
\wedge \tilde{f}(\gf_\ga(k_0)) = \tilde{c}(\ga)$ since $k_0 > m_{\ga_0}$. \qed

\section{Finitely many primes}
\label{4monu}

The proof of Theorem~\ref{1monuthm} uses infinitely many primes.
Otherwise said, the type of the (torsion-free rank one) non-free
quotients $G_{\gd + 1}/G_\gd$ in that construction is $(1, 1, 1,
\ldots )$. We may ask what happens if we are allowed only finitely
many primes. For example, we may consider modules over $\oZ_{(P)}$
(where $P$
 is a set of primes and $\oZ_{(P)}$ denotes the rationals whose
denominators in reduced form are not divisible by an element of $P$)
and ask whether the main theorem, Theorem~\ref{main},
 holds.  If $P$ is infinite, i.e.,
$\oZ_{(P)}$ has infinitely many primes, then our proofs apply and
there is a non-free hereditarily separable $\oZ_{(P)}$-module of
cardinality $\ha_1$ if and only if there is a ladder system on $\goo$
which satisfies monochrome uniformization for $\go$ colours.  On the other hand if the
cardinality of $P$ is finite but at least two,
we can  show that Theorem~\ref{main} does not
hold, and Case 4 in the Introduction is impossible. In fact, this section is
devoted to proving the following result:

\begin{theorem}
\label{4monuthm}
Suppose $R$ is a countable p.i.d. with only finitely many but at least
$2$ primes. If
there is an hereditarily separable $R$-module of cardinality $\hao$
which is not free, then there is a Whitehead $R$-module
of cardinality $\hao$  which is not
free.
\end{theorem}

\proof
The method of proof is to show that if there is an hereditarily
separable $R$-module of cardinality $\hao$ which is not free, then there
is a ladder system on a stationary subset of $\goo$ which satisfies
$2$-uniformization.
We first prove
 that
 \labpar{\dagger}{there is a tree-like ladder system $\gh =
\set{\gh_\gd}{\gd \in S}$ on a stationary subset $S$ of $\lim(\goo)$
such that for every $2$-colouring $c = \set{c_\gd}{\gd \in S}$ of $\gh$,
there is a function $f \colon \goo \times \go \to 2$ such that for all
$\gd \in S$ there exists $m_\gd \in \go$ such that $f(\gh_\gd(n), m_\gd)
= c_\gd(n)$ for all $n \in \go$.}

\noindent
 Let $N$ be an hereditarily separable
$R$-module of cardinality $\hao$. As in the proof of
Theorem~\ref{2monuthm}, we  write $N = \cup _{\alpha < \omega _1}
N_\alpha $ as a union of a continuous chain of countable free pure
submodules where there is a stationary subset $S$ of $\omega _1$,
consisting of limit ordinals, and an integer $r \geq 0$ such that for
all $\delta \in S$, $N_{\delta +1}/N_\delta $ is non-free of rank $r +
1$ and every subgroup of $N_{\delta +1}/N_\delta $ of rank $r$ is free.
There is a pure free subgroup $F_\delta /N_\delta $ of $N_{\delta
+1}/N_\delta $ of rank $r$ such that $N_{\delta +1}/F_\delta $ is rank 1
and non-free.

It follows from the fact that there are only finitely many primes that the type
of $N_{\delta +1}/F_\delta $ is
 $(t_1, t_2, \ldots,
t_n)$ where each $t_i$ is either $0$ or $\infty$
and at least one $t_i =
\infty$. Thus without loss of generality we may assume that there is a
fixed prime $p \in R$ such that for all $\gd \in S$ there exists
 $\{y^\delta _\ell \colon 0 \leq \ell \leq r\} \se N_{\gd + 1}$
 such that $\{y^\delta _\ell + N_\delta \colon 0 \leq \ell \leq r -1\}$
 is a basis of $F_\delta/N_\gd $, $y^\delta _r \notin F_\gd$ and
 $y^\delta _r + F_\delta $ is  $p$-divisible in $N_{\delta + 1}/F_\gd$.
 Then $N_{\delta +1}$ contains elements
$z^\gd_n$ ($n \in
 \omega$), where $z^\gd_0 = y^\gd_r$ and the $z^\gd_n$ ($n \geq 1$) satisfy equations
\labpar{\blacktriangle_n}{$p z^\gd_n
 = \sum _{\ell \leq r} s^{\delta ,n}_\ell y^\delta _\ell + g^\gd_n +
 \sum_{j < n}r^{\gd, n}_j z^\gd_j $}

\noindent
where $g^\gd_n \in  N_\delta $, $r^{\delta, n}_\ell$, $s^{\delta ,n}_\ell
\in  R$, and no element of  $N_{\gd + 1}/\langle  F_\gd \cup
\set{z^\gd_n}{n \in \go} \rangle$ has order $p$.

Define functions $\varphi _\delta $ on $\omega $ for each $\delta  \in  S$ by:
$$
\varphi _\delta (n) = \langle g^\gd_m,
s^{\delta ,m}_\ell ,  r^{\delta, m}_j \colon \ell  \leq  r,  m \leq
n, j < m\rangle .
$$
Let $c$ be a $2$-colouring of $\Phi
= \{\varphi _\delta \colon \delta \in S\}$. Following the pattern of the
proof of Theorem~\ref{2monuthm},  we will use $c$ to define a
subgroup $B$ of $N$. We begin  by letting $\{x_n \colon n \in \omega \}$
be a basis of $N_0$, and letting $B_0$ be the subgroup of $N_0$ generated
by $\{px_0\} \cup \{px_{n+1} - x_n\colon n \in \omega \}$.  Also, let
$A = \{t_n\colon n \in \omega \} \subseteq N_0$ be a complete set
of representatives of $N_0/B_0$ such that $t_0 = 0$ and
for each
 $a \in A$, fix
an element $[p, a] \in \go$ such that
 $p t_{[p,a]} + B_0 = a + B_0$.

Assume we have defined $B_\gd$ so that $B_\gd + N_0 = N_\gd$ and for all
$\gb < \gd$, $B_\gd \cap N_\gb = B_\gb$. We now define
$B_{\gd, m}$ inductively so that $z^\gd_{m} \in B_{\gd, m} + N_0$. Let
$B_{\gd, 0}$ be generated by $B_\gd \cup \{y^\gd_0, \ldots, y^\gd_r\}$.
If $B_{\gd, m-1}$ has been defined, we have $pz^\gd_m \in B_{\gd, m-1} +
N_0$, so $pz^\gd_m \equiv a^\gd_m \pmod{B_{\gd, m-1}}$ for some
$a^\gd_m \in A$. Let $$B_{\gd, m} = B_{\gd, m-1} + R(z^\gd_m - t_{[p,
a^\gd_m]} - c_\gd(m)x_0).$$ Having defined $B_{\gd, m}$ for all $m$, we
can extend $\bigcup_{m \in \go} B_{\gd, m}$ to $B_{\gd + 1}$ such that
$B_{\gd + 1} + N_0 = N_{\gd + 1}$ and $B_{\gd + 1} \cap N_\gd = B_\gd$.

Finally, let $B = \bigcup_{\ga < \goo}B_\ga$ and fix a projection $h\colon B
\rightarrow Rpx_0$ and  a well-ordering, $\prec$, of $R^{r+1} $ of order
type $\omega $. Extend $h$ to a homomorphism, also denoted $h$, from $N$
into $Qpx_0$, where $Q$ is the quotient field of $R$. Given $\nu$ of the
form $\gf_\gd(n)$ and $m \in \go$, we are going to define $f(\nu, m)$.
Let $\langle w^m_\ell \colon \ell \leq r \rangle$ be the $m$th tuple in
$R^{r+1}$ according to $\prec$. We shall suppose that
\labpar{\#_m}{$h(y^\gd_\ell) =
w_\ell^m px_0$}
 for $\ell \leq r$ for some $\gd$ such that $\nu = \gf_\gd(n)$,
and show that under this supposition (and with the information given by
$\nu$) we can compute $c_\gd(n)$; we will then define this value of
$c_\gd(n)$ to be $f(\nu, m)$. Since one of our suppositions ($\#_m$)
 about the values of $h(y^\gd_\ell)$ must be right, ($\dagger$)
will be proved.

The proof is by induction on $k \leq n$ that we can compute
$h(z^\gd_k)$, $a^\gd_k$, and $c_\gd(k)$.
In fact, for $0 < k \leq n$ we have an equation
\labpar{\blacktriangle_k}{$p z^\gd_k
 = \sum _{\ell \leq r} s^{\delta ,k}_\ell y^\delta _\ell + g^\gd_k +
 \sum_{j < k}r^{\gd, k}_j z^\gd_j $}
 satisfied by $z^\gd_k$. Since by
 induction and our supposition we know  the value of $h$ for
 all the elements on the right-hand side, we can compute $h(z^\gd_k)$.
 Since by induction we also know $c_\gd\rest k$, we know $B_{\gd, k-1}$,
 so we can calculate $a^\gd_k$ ($\equiv pz^\gd_k \pmod{B_{\gd, k-1}}$).
 Finally, we know that
  $$h(z^\gd_k - t_{[p, a^\gd_k]} - c_\gd(k)x_0) =
 h(z^\gd_k - t_{[p, a^\gd_k]}) - c_\gd(k)x_0 $$
  belongs to $Rpx_0$, and
 we know $h(z^\gd_k - t_{[p, a^\gd_k]})$ by induction (because we know
 $h\rest N_0$).
  Now $h(z^\gd_k - t_{[p, a^\gd_k]}) - x_0$ and $h(z^\gd_k - t_{[p,
  a^\gd_k]})$
  cannot both belong to $Rpx_0$.
If $h(z^\gd_k - t_{[p, a^\gd_k]})$ belongs to
 $Rpx_0$, $c_\gd(k)$ must equal $0$;  otherwise let $c_\gd(k) = 1$. (If
 the latter value does not make $h(z^\gd_k - t_{[p, a^\gd_k]} -
 c_\gd(k)x_0)$ belong to $Rpx_0$, then our supposition must have been
 wrong, and we can let $f(\nu, m)$ be arbitrary.)

This completes the proof of ($\dagger$). At this point we use the
assumption that there
 are at least two primes.
Then  the proof of
necessity, i.e. of Theorem~\ref{2monuthm}, is still valid.
(Referring to the last paragraph of that proof, we use the fact that
there are two primes when we assert that $H/K$ is not divisible.)
Moreover, there is a single ladder system $\gh$ which satisfies the
property of ($\dagger$) as well as monochromatic uniformization for $\go$
colours. (Indeed, by reducing to a smaller set we can assume that the
same set $S$ is used in both the proof of ($\dagger$) and the proof of
Theorem~\ref{2monuthm}; then we can let $\gh$ be a ladder system derived
from functions $\gf_\gd$ which give combined information about the
equations used in the proof of ($\dagger$) and the equations used in the
proof of Theorem~\ref{2monuthm}.)

 Given a $2$-colouring $c$ of $\gh$, let $f$ be as in ($\dagger$).
Define a monochromatic colouring $c'$ of $\gh$ by: $c'(\gd) = m_\gd$
where $m_\gd$ is such that $f(\gh_\gd(n), m_\gd) = c_\gd(n)$ for all
$n \in \go$. There is a uniformization $\langle g, g^*
\rangle$ of $c'$. Define $h \colon \goo \to 2$ by: $h(\nu) = f(\nu,
g(\nu))$.  Then for all $ \gd \in S$
for sufficiently large $n$, $$h(\gh_\gd(n)) = f(\gh_\gd(n),
g(\gh_\gd(n))) = f(\gh_\gd(n), m_\gd) = c_\gd(n).$$ \qed

The third author, in \cite[Thm. 3.6]{Sh105}, claimed to prove  that
if the non-freeness of $G$ involves only finitely many primes, then $G$
is hereditarily separable if and only if $G$ is Whitehead. However, the
proof given seems to be irredeemably defective. We do not know if the
result claimed is true. Thus we  still have
  the following open questions:

\begin{quote}
If $R$ is a countable p.i.d. with exactly one prime, does
Theorem~\ref{4monuthm} hold?
If $R$ has finitely many primes, is every hereditarily separable $R$-module of
cardinality $\hao$ a Whitehead module? If not, find a combinatorial
equivalent, analogous to Theorem~\ref{main}, to the existence of a
hereritarily-separable $R$-module which is not a Whitehead module.
 \end{quote}

\end{document}